\documentclass[12pt]{amsart}
\usepackage{amssymb,amsthm} 
\usepackage{pxfonts}   

\newtheorem{thm}{Theorem}

\newtheorem{lemma}{Lemma}
\theoremstyle{definition}
\newtheorem{rem}{Remark}

\newtheorem{claim}{Claim}
\theoremstyle{remark}
\newtheorem*{ack}{Acknowledgements}

\newcommand{\comm}[1]{}

\def\R{\mathbb{R}}

\def\N{\mathbb{N}}
\def\eps{\varepsilon}
\def\cB{\mathcal{B}}
\def\cC{\mathcal{C}}

\def\cM{\mathcal{M}}

\def\cR{\mathcal{R}}

\def\cV{\mathcal{V}}

\def\id{{\operatorname{id}}}
\renewcommand{\epsilon}{\varepsilon}
\renewcommand{\setminus}{\smallsetminus}
\renewcommand{\emptyset}{\varnothing}

\begin{document}

\title[A generic map has no acim]{A generic $C^1$ map has no absolutely continuous invariant probability measure}

\author[A.~Avila]{Artur Avila}
\address{CNRS UMR 7599,
Laboratoire de Probabilit\'es et Mod\`eles al\'eatoires.
Universit\'e Pierre et Marie Curie--Bo\^\i te courrier 188.
75252--Paris Cedex 05, France}
\urladdr{www.proba.jussieu.fr/pageperso/artur/}
\email{artur@ccr.jussieu.fr}

\author[J.~Bochi]{Jairo Bochi}
\address{Instituto de Matem\'atica -– UFRGS -- Porto Alegre, Brazil}
\urladdr{www.mat.ufrgs.br/{\textasciitilde}jairo}
\email{jairo@mat.ufrgs.br}

\date{May 29, 2006}

\maketitle

\section{Statement}

Let $M$ be a smooth compact manifold
(maybe with boundary, maybe disconnected) of any dimension $d \ge 1$.
Let $m$ be some (smooth) volume probability measure in $M$.
Let $C^1(M,M)$ be the set of $C^1$ maps $M \to M$, endowed with the $C^1$ topology.
Given $f\in C^1(M,M)$, we say that $\mu$ is an \emph{acim for $f$}
if $\mu$ is an $f$-invariant probability measure which is absolutely continuous
with respect to~$m$.

\begin{thm}\label{t.main}
The set $\cR$ of $C^1$ maps $f:M \to M$
which have no acim is a residual (dense $G_\delta$) subset of $C^1(M,M)$.
\end{thm}

Since the set the set of all expanding maps
and the set of all diffeomorphisms
are open subsets of $C^1(M,M)$, we have the following immediate consequences:
\begin{enumerate}
\item The $C^1$-generic expanding map has no acim.
\item The $C^1$-generic diffeomorphism has no acim.
\end{enumerate}

Result (i) was previously obtained in the case $M$ is the circle by Quas~\cite{Quas}.
Of course, (i) does not hold in the $C^{1+\text{H\"older}}$ topology.

It seems possible that result (ii) holds in higher topologies.
An old result by Livsic and Sinai implies that the $C^\infty$-generic
Anosov map has no acim, see~\cite{LivsicSinai}, also~\cite{Chernov}.
(In fact, the existence of a single periodic point of the Anosov map over which
the Jacobian is different from~$1$ prohibits the existence of an acim.)
On the other hand, the existence of acim is certainly not rare
(in the probabilistic sense) among smooth enough diffeomorphisms of tori
close to translations (by KAM).

In the course of the proof, we will need a generalization of the usual
Rokhlin tower lemma to non-invariant measures.
That result, theorem~\ref{t.rokhlin}, may be of independent interest.

\section{Proof}

In \S\S\ref{ss.criterium}--\ref{ss.slicing} we give some results
that are added up to prove theorem~\ref{t.main} in \S\ref{ss.proof}.

\subsection{Criterium for existence of acim}\label{ss.criterium}

The following shows that the set of
maps that do not have an acim is a $G_\delta$ subset of $C^1(M,M)$.

\begin{lemma}\label{l.no acim criterium}
A map $f\in C^1(M,M)$ has no acim iff for every $\eps>0$
there exists a compact set $K \subset M$ and $N \in \N$ such that
$$
m(K) > 1 - \eps \quad \text{and}\quad m(f^N(K)) < \eps.
$$
\end{lemma}

\begin{proof}
Assume that $f$ has an acim $\mu$.
Let $\eps>0$ be such that $m(Z) \le \eps$ implies $\mu(Z) < 1/2$.
Now assume that $K \subset M$ is a compact set
such that $m(f^N K) < \eps$ for some $N \in \N$.
Then $\mu(K) \le \mu(f^N K) < 1/2$, so
$m \left(M \setminus K\right) > \eps$.

Next assume that $f$ has no acim.
Let $\mu$ be a limit point of the sequence of measures
$\frac{1}{n} \left( m + f_* m + \cdots + f^{n-1}_* m\right)$;
then $\mu$ is $f$-invariant.
Let $\mu = \mu_\text{ac} + \mu_\text{sing}$ be the Lebesgue decomposition
of $\mu$ relative to $m$.
Since $f$ is $C^1$, $f_* \mu_{\text {sing}}$ is singular, and it follows
that $\mu_\text{ac}$ and $\mu_{\text {sing}}$ are $f$-invariant.
But $f$ is assumed to have no acim, so $\mu = \mu_\text{sing}$.
Thus there exists $Z \subset M$ such that $m(Z) = 1$ and $\mu(Z) = 0$.
Given any $\eps>0$,
take a compact set $L$ and an open set $V$ such that $L \subset Z \subset V$,
$m(L) > 1 - \eps$ and $\mu(V) < \eps/2$.
Let $\phi$ be a continuous function such that $\chi_L \le \phi \le \chi_V$.
For some sequence $n_j \to \infty$ we have
$$
\frac{1}{n_j} \sum_{i=0}^{n_j-1} \int f^i \circ \phi \; dm \to \int \phi \; d\mu < \eps/2 \, .
$$
In particular, there exists $N$ such that
$m(f^{-N} L) \le \int f^N \circ \phi \; dm < \eps/2$.
Take a compact $K \subset M \setminus f^{-N}L$ such that
$m(M\setminus K)<\eps$.
Then $m(f^N K) \le m(M\setminus L) < \eps$.
\end{proof}

\subsection{A non-invariant Rokhlin lemma}\label{ss.rokhlin}

\begin{thm}\label{t.rokhlin}
Let $f:M \to M$ be a $C^1$ endomorphism of a compact manifold,
and let $m$ be normalized Lebesgue measure.
Assume that $m(C_f \cup P_f)=0$,
where $C_f$ is the set of critical points and $P_f$ is the set of periodic points.
Given any $\epsilon_0>0$ and $n_0$, $\ell \in \N$ with $\ell \le n_0$,
there exists a measurable set $U \subset M$ such that
$f^{-i}(U) \cap U = \emptyset$ for $1 \leq i < n_0$,
\begin{equation}\label{e.rokhlin conclusion}
\sum_{i=0}^{n_0-1} m \left( f^{-i}(U) \right) > 1 -\epsilon_0,
\quad \text{and} \quad
\sum_{i=0}^{\ell-1} m \left( f^{-i}(U) \right) < \frac{\ell}{n_0} + \eps_0 \,.
\end{equation}
\end{thm}

Notice that if the map $f$ were assumed to preserve the measure $m$,
the theorem would be an immediate consequence
of the well-known Rokhlin lemma (for non-invertible maps, see \cite{rokhlin}).

The proof of theorem~\ref{t.rokhlin} will occupy the rest of this subsection.
Let $f$ be fixed from now on.

Let $\cM$ be the $\sigma$-algebra of measurable sets.
Since $f$ is $C^1$, $Y \in \cM$ implies $f(Y) \in \cM$.

Given $Z \in \cM$, we denote
$$
\widehat{Z}=\bigcup_{i=0}^\infty f^{-i}(Z) \, .
$$
We say that $Z \in \cM$ is $N$-good (where $N\in\N$) if $Z \cap f^{-i}(Z)=\emptyset$ for $0 \leq i<N$.

\begin{claim}\label{cl.union}
If $A$ and $B$ are $N$-good sets then the set
$$
C = \left(A \setminus \widehat B \right) \cup
\left(B \setminus \left(A \setminus \widehat B \right)^{\mathord{\wedge}}\right)
$$
is $N$-good and satisfies $\widehat C \supset A \cup B$.
\end{claim}


\begin{proof}
Let $A'= A \setminus \widehat B$;
then $(A' \cup B)^{\mathord{\wedge}} = (A \cup B)^{\mathord{\wedge}}$.
Let $B'= B \setminus \widehat{A'}$;
then $(A' \cup B')^{\mathord{\wedge}} = (A' \cup B)^{\mathord{\wedge}}$.
That is, the set $C = A' \cup B'$ satisfies $\widehat C = (A \cup B)^{\mathord{\wedge}}$.
Using that $A'$ and $B'$ are $N$-good, $A' \cap \widehat{B'} = \emptyset$,
and $\widehat{A'} \cap B' = \emptyset$, we see that $C$ is $N$-good.
\end{proof}

We say that $Z\in \cM$ is $N$-saturated if $f^{-N}(f^N(Z)) = Z$.
The $N$-saturated sets form the $\sigma$-algebra $f^{-N}\cM$.

\begin{claim}\label{cl.cover}
For each $N \in \N$ there exists a countable cover
(modulo sets of zero $m$-measure)
$M = \bigcup B_k$ such that each $B_k$ is $N$-good and $N$-saturated.
\end{claim}

\begin{proof}
Since $m(P_f)=0$, there is a countable cover $M = \bigcup A_k$,
where the sets $A_k$ are $N$-good.
Take $B_k = f^{-N}(A_k)$.
\end{proof}

\begin{claim}\label{cl.W}
For every $\eps>0$ and $N\in\N$ there exists a set $W$
which is $N$-good, $N$-saturated, and $m(\widehat W)>1-\eps$.
\end{claim}

\begin{proof}
Let $B_k$ be the sets given by claim~\ref{cl.cover}.
Define inductively sets $C_k$: take $C_1 = B_1$,
and for $k>0$, let $C_{k+1}$ be the $N$-good set given by claim~\ref{cl.union}
such that $\widehat{C_{k+1}} \supset C_k \cup B_{k+1}$.
Then for all $k$ we have that $C_k$ is $N$-saturated and
$\widehat{C_k} \supset \bigcup_{j=1}^k B_j$.
Finally, take $W=C_k$ for some large~$k$.
\end{proof}

\begin{claim}\label{cl.V}
For every $\eps>0$ and $N \in \N$, there exists a $N$-good set $V$ such that
$m(V)<\epsilon$ and $m(\widehat V)>1-\epsilon$.
\end{claim}

\begin{proof}
Increasing $N$ if necessary, we assume $N> 1/\epsilon$.
Take $W$ as in claim~\ref{cl.W}.
Notice that the sets $W$, $f(W)$, \ldots, $f^{N-1}(W)$ are disjoint.
Take $0 \leq i \leq N-1$ such that $m(f^i(W)) \le 1/N$.
Let $V = f^i(W)$; then $\widehat V \supset \widehat W$.
Since $W$ is $N$-good and $N$-saturated, $V$ is $N$-good.
\end{proof}

\begin{claim}\label{cl.abs cont}
For any $i \ge 0$, $f^i_* m$ is absolutely continuous with respect to~$m$.
\end{claim}

\begin{proof}
Clearly it suffices to consider $i=1$.
Let $Z \in \cM$ be such that $m(f^{-1}(Z))>0$.
Since $m(C_f)=0$, we can find an open set $U \subset M \setminus C_f$
such that $f|U$ is a $C^1$-diffeomorphism and $m(f^{-1}(Z) \cap U)>0$.
Then $f(f^{-1}(Z) \cap U)$, and hence $Z$, has positive measure.
\end{proof}

\begin{proof}[Proof of theorem~\ref{t.rokhlin}]
Let $\ell$, $n_0$, and $\eps_0$ be given.
By claim~\ref{cl.abs cont}, there exists $\eps>0$ be such that
$$
Z \in \cM, \ m(Z)<\eps \  \Rightarrow \
m \left( \bigcup_{i=0}^{2n_0-1} f^{-i} Z \right) < \frac{\epsilon_0}{2} \, .
$$
Let $V$ be given by claim~\ref{cl.V} with $N=n_0$.
For $i\ge 0$, let $V_i = f^{-i}(V)$ and
$$
V_i^* = V_i \setminus \bigcup_{j=0}^{i-1} V_j \, .
$$
For each $0 \le j < n_0$, let
$$
S_j = \sum_{k=j}^{j+\ell-1} \sum_{\substack{i \ge 0 \\ i=k \bmod n_0}} m(V_i^*) \, .
$$
We have
$\sum_{j=0}^{n_0-1} S_j = \ell \cdot m(\widehat V)$,
so there exists some $j_0$ for which $S_{j_0} \le \ell/n_0$.
Define
$$
U=\bigsqcup_{\substack{i \ge n_0 \\ i=j_0 \bmod n_0}} V_i^* \, .
$$
Noticing that
$f^{-j}(V_i^*) \subset V_{i+j}^* \cup V_0 \cup V_1 \cup \ldots \cup V_j$
for $0 \le j < n_0$,
we see that $U$ is $n_0$-good.
Also,
$$
m \left( \bigsqcup_{j=0}^{\ell-1} f^{-j} (U) \right) \le
m(V_0 \cup \cdots \cup V_{\ell-1}) + S_{j_0} <
\frac{\ell}{n_0} + \frac{\eps_{0}}{2} \, .
$$
Finally, since $f^{-j}(V_i^*) \supset V_{i+j}^*${\,}, we have
$$
m \left( \bigsqcup_{j=0}^{n_0-1} f^{-j} (U) \right)
\ge m \left( \bigsqcup_{i=2n_0}^{\infty} V_i^* \right)
> 1 - \eps - \frac{\eps_{0}}{2} \ge 1 - \eps_0 \, .
$$
\end{proof}

\begin{rem}
We used only the following assumptions about $f$ and $m$:
\begin{itemize}
\item $(M,\mathcal{M},m)$ is a Lebesgue space and $f:M \to M$ is measurable;
\item $f$ is aperiodic: $m(P_f)=0$;
\item $f$ is non-singular with respect to $m$: for $Y \in \mathcal{M}$, 
we have $m(Y)=0$ if and only if $m(f^{-1}(Y))=0$;
\item $f$ is forward-measurable: $Y \in \mathcal{M}$ implies 
$f(Y) \in \mathcal{M}$. 
(In fact, we can always replace $f$ by a isomorphic copy 
which is forward-measurable: see \cite{Rohlin}.)
\end{itemize}\end{rem}

\begin{rem}[Addendum to theorem~\ref{t.rokhlin}]\label{r.open}
The set $U$ can be taken open,
and with $f^{-i}(\overline U) \cap \overline U = \emptyset$, $0\le i < n_0$.

Indeed, take a compact set $K \subset U$
with $m(U \setminus K)$ very small.
Then take an open set $U_0 \supset K$
with $m(U_0 \setminus K)$ very small and such that
$\overline{U_0}$, \ldots, $f^{-n_0+1}(\overline{U_0})$ are disjoint.
Bearing in mind claim~\ref{cl.abs cont},
we see that~\eqref{e.rokhlin conclusion} holds with $U_0$ in the place of $U$.
\end{rem}

\subsection{Linearization}\label{ss.linearization}

Fix an atlas of $M$ formed by charts that take the
restricted volume on $M$ to Lebesgue measure on $\R^d$.  Fix also a family
of pairs $(A_i,\phi_i)$ such that the $A_i \subset M \setminus \partial M$
are disjoint open sets compactly
contained in the domain of the chart $\phi_i$, and $\sum m(A_i)=1$.  We call
the $A_i$ {\it basic blocks}.

We shall say that a map $f:M \to M$
is \emph{locally linear on an open set $V \subset M \setminus \partial M$}
if for each connected component $W$ of $V$, there exists
both $W$ and $f(W)$ are contained in basic blocks
and if under the corresponding change of coordinates the map
$f:W\to f(W)$ becomes the restriction
of an affine map $\R^d \to \R^d$.

\begin{lemma}\label{l.linearization}
If $f:M \to M$ is a $C^1$ map and $U \subset M \setminus \partial M$ is open then
for every $\gamma>0$
there exists a $C^1$-map $\tilde f: M \to M$
which is $C^1$-close to $f$ and equals $f$ outside $U$, 
and there exists an open set $V \subset U$ such that
$m(V)/m(U) > 1 - \gamma$ and $\tilde f$ is locally linear on $V$.
Furthermore, if the set $C_f$ of critical points of $f$ has zero Lebesgue measure
then $\tilde f$ can be taken to be a local diffeomorphism on $V$.
\end{lemma}

\begin{proof}
Up to reducing $U$ a little, we can assume each connected component of $U$,
as well as its image by $f$, is contained in a basic block.

To simplify writing, from now on we assume $U\subset \R^d$, and $m$ is Lebesgue measure on $\R^d$.
Given $\gamma>0$, let $\delta>0$ be such that $(1-\delta)^d < \gamma/2$.
Fix a $C^1$ bump function $\rho: \R^d \to [0,1]$ such that
$\rho(x)=1$ if $\|x\| \le 1-\delta$,
$\rho(x)=0$ if $\|x\| \ge 1$
(where $\|\mathord{\cdot}\|$ is the euclidian norm on $\R^d$).

Let $r_0>0$ be small.
By Vitali's lemma, we can find finitely many disjoint balls
$B(p_i,r_i) \Subset U$ of radii $r_i<r_0$, such that the Lebesgue measure
of their union is greater than $(1-\gamma/2)m(U)$.
Define $\tilde f$ on each $B(p_i,r_i)$ by
$$
\tilde f(x) =
f(x) + \rho(r_i^{-1}(x-p_i)) \cdot [-f(x) + f(p_i) + Df(p_i) \cdot (x-p_i)] \, ,
$$
and $\tilde f=f$ on $U \setminus \bigsqcup_i B(p_i,r_i)$.

Then $\tilde f$ is locally linear on \mbox{$V = \bigsqcup_i B(p_i, (1-\delta) r_i)$}.
If $r_0$ is sufficiently small, then $\tilde f$ is $C^1$-close to~$f$.

If $m(C_f)=0$ we take each $p_i$ such that $Df(p_i)$ is an isomorphism.
\end{proof}

\subsection{Perturbation of a sequence of linear maps} \label{ss.slicing}

\begin{lemma}\label{l.slicing}
Given $\eps>0$ and $0<\delta<1$, there exists $k \in \N$
such that given any of sequence linear isomorphisms
$$
\R^d \xrightarrow{L_n} \R^d \xrightarrow{L_{n-1}} \cdots \xrightarrow{L_1} \R^d \, ,
$$
with $n \ge k$, there exists $\tau_0>0$ such that
for any $0<\tau<\tau_0$ the following holds true:

Define boxes
\begin{alignat*}{2}
U_0 &= [-1,1]^{d-1} &&\times [-\tau,\tau], \\
V_0 &= [-(1-\delta), 1-\delta]^{d-1} &&\times [-(1-\delta)\tau, (1-\delta)\tau], \\
W_0 &= [-1, 1]^{d-1} &&\times [-\delta\tau, \delta \tau].
\end{alignat*}
Define also $U_i = L_i^{-1} U_{i-1}$, $V_i = L_i^{-1} V_{i-1}$, $W_i = L_i^{-1} W_{i-1}$,
for $1 \le i \le n$.
Then there exist $C^1$-diffeomorphisms $H_i:\R^d \to \R^d$
with derivative $\eps$-close to $\id$,
and with $H_i = \id$ outside $U_i$,
such that for all $i$ with $k \le i \le n$ we have
\begin{equation}\label{e.shrink}
L_{i-k+1} \circ H_{i-k+1} \circ \cdots \circ L_{i-1} \circ H_{i-1} \circ L_i \circ H_i (V_i) \subset W_{i-k} \, .
\end{equation}
Moreover, for $1 \leq i \leq n$, $H_i$ only depends on $\epsilon$, $\delta$,
$\tau$, and $L_1,...,L_i$ (but not on $L_{i+1},...,L_n$).
\end{lemma}

In the following proof of lemma~\ref{l.slicing}, we will assume $d \ge 2$,
leaving for the reader the easy adaptation to the case $d=1$
(where any $\tau_0$ works).

In the proof we will need
lemmas~\ref{l.basic perturb} and \ref{l.trunca} below.
We write $\R^{d-1} = \R^{d-1} \times \{0\} \subset \R^d$.
Also, we call a subset $B$ of a finite-dimensional vector space $V$ a \emph{ball}
if there exists a norm on $V$ such that $B$ is the closed unit ball on $V$ with respect to that norm.

\begin{lemma}\label{l.basic perturb}
For every $\eps>0$ and $0<\delta<1$, there exists $0<\kappa<1$ with the following properties:
Given any ball $\cC \subset \R^{d-1}$, there exists
$\tau^*>0$ such that
if $0<\tau<\tau^*$ then there exist a diffeomorphism $H:\R^d \to \R^d$ satisfying the following:
\begin{itemize}
\item $H$ has derivative $\eps$-close to the identity;
\item $H$ equals the identity outside $\cC \times [-\tau,\tau]$;
\item if $(z,t) \in (1-\delta)(\cC \times [-\tau,\tau])$ then $H(z,t) = (z, \kappa t)$.
\end{itemize}
\end{lemma}

\begin{proof}
Given $\eps$ and $\delta$, let $\kappa$ be such that
$$
(1-\kappa) (1 + 2\delta^{-1})<\eps \, .
$$
Now let $\cC = \{z \in \R^{d-1}; \; \|z\|_* \le 1\}$,
where $\|\mathord{\cdot}\|_*$ is a norm in $\R^{d-1}$.
Let $C>0$ be such that $\|v\|_* \le C \|v\|$, where
$\|\mathord{\cdot}\|$ is euclidian norm.
Let $\tau^* = C^{-1}$.

Then, given $0<\tau<\tau^*$,
take a bump function $\rho:\R \to [0,1]$ such that
$\rho(x)=1$ for $|x|\le 1-\delta$,
$\rho(x)=0$ for $|x|\ge 1$,
and $|\rho'| < 2 \delta^{-1}$.
Define
$$
H(z,t) = \left(z, \; \left[1-(1-\kappa) \cdot \rho(\tau^{-1} t) \cdot \rho(\|z\|_*) \right] t \right) \, ,
\quad z \in \R^{d-1}, \ t\in \R.
$$
Then
$$
\left| \frac{\partial H}{\partial t}  - 1 \right| \le
(1-\kappa) \cdot \rho(\|z\|_*) \cdot
\left[ \rho(\tau^{-1} t)  +
\left| \rho'(\tau^{-1} t) \cdot \tau^{-1} t \right| \right]
< \eps \, .
$$
And if $v\in \R^{d-1}$ then
\begin{align*}
\left\| D H(z,t) \cdot (v,0) \right\| &\le
|t| \cdot (1-\kappa) \cdot \rho(\tau^{-1} t) \cdot \big| \rho'(\|z\|_*) \big| \cdot \|
v\|_*\\
\nonumber
&\leq \tau^* (1-\kappa) \cdot 2 \delta^{-1} \cdot C \|v\| \le \eps \|v\| \, .
\end{align*}
\end{proof}

Let $e_d = (0,\ldots,0,1) \in \R^d$ and $\langle \mathord{\cdot}, \mathord{\cdot} \rangle$
be the standard inner product on $\R^d$.
The easy proof of the following lemma is left to the reader:

\begin{lemma} \label{l.trunca}
Let $L:\R^d \to \R^d$ be a linear isomorphism
such that $L(\R^{d-1})=\R^{d-1}$.
Let $\beta = |\langle L(e_d), e_d \rangle|$.
Let $\cC \subset \R^{d-1}$ be a ball, and let $\lambda>1$.
Then there exists $\tau'>0$ such that
for any $0 < \tau < \tau'$ we have:
$$
(\lambda^{-1} L(\cC)) \times [-\beta \tau, \beta \tau]
\subset L\big( \cC \times [-\tau,\tau] \big) \subset
(\lambda L(\cC)) \times [-\beta \tau, \beta \tau] \, .
$$
\end{lemma}

\begin{proof}[Proof of lemma~\ref{l.slicing}]
Let $\kappa = \kappa(\eps,\delta/2)>0$ be given by lemma~\ref{l.basic perturb}.
We take $k \in \N$ such that $\kappa^k < \delta/(1-\delta)$.

Now take $n \ge k$ and $L_1$, \ldots, $L_n$ as in the statement of the lemma.
Rotating coordinates if necessary, we can assume
$L_i \cdot \R^{d-1} = \R^{d-1}$ for all $i$.
Let $\cC_0 = [-1,1]^{d-1}$ and
$\cC_i = L_i^{-1} \cdots L_1^{-1} \cdot \cC_0$ for $i \ge 1$.
Let $\alpha_0 = 1$ and $\alpha_i = \left| \langle L_i^{-1} \cdots L_1^{-1} e_d, e_d \rangle \right|$.
Write, for $a>0$ and $b>0$,
$$
\cB_i[a, b] = (a \cC_i) \times [-\alpha_i b, \alpha_i b] \, .
$$
Fix $\lambda>1$ so that
\begin{equation}\label{e.lambda}
\lambda^{3n} < (1-\delta)^{-1} (1- \delta/2) \, .
\end{equation}

Let $\tau_i^*>0$ be associated to the ball $\cC = \cC_i$ by lemma~\ref{l.basic perturb}.
Let $\tau'_i>0$ be associated to the linear map $L_i^{-1}$,
the ball $\cC_i$, and $\lambda$ by lemma~\ref{l.trunca}.
Let
$$
\tau_0 = \lambda^{-2n} \min \left\{\alpha_i^{-1} \tau^*_i, \ \alpha_i^{-1} \tau'_i \, ; \;  1\le i \le n\right\}.
$$

By lemma~\ref{l.trunca},
\begin{equation}\label{e.chatice}
\cB_i [\lambda^{-1} a , b] \subset L_i^{-1} \left( \cB_{i-1}[a,b] \right) \subset
\cB_i [\lambda a, b] , \quad \text{provided }\frac{b}{a}< \frac{\tau'_i}{\alpha_i} \, .
\end{equation}
Therefore
\begin{equation}\label{e.sono}
\frac{b}{a}< \lambda^n \tau_0 \ \Rightarrow \
\cB_i [\lambda^{-n} a, b] \subset
(L_i \circ \cdots \circ L_1)^{-1} \left( \cB_0[a, b] \right) \subset
\cB_i [\lambda^n a, b] \, .
\end{equation}

Now let $0<\tau<\tau_0$ be fixed,
and let $U_0$, $V_0$, $W_0$ be as in the statement of the lemma, that is,
$$
U_0 = \cB_0[1,\tau], \quad
V_0 = \cB_0[(1-\delta), (1-\delta)\tau], \quad
W_0 = \cB_0[1, \delta \tau].
$$

Let $\tilde H_i:\R^d \to \R^d$ be the diffeomorphism
supported on
$\cB_i(1, \lambda^n \tau) = \cC_i \times [-\alpha_i \lambda^n \tau , \alpha_i \lambda^n \tau]$
given by lemma~\ref{l.basic perturb}.
We define $H_i$ by
$$
H_i(x) = \lambda^{-n} \cdot \tilde H_i (\lambda^n x)  \, .
$$
Then $DH_i$ is $\eps$-close to $\id$.
Also, $H_i$ equals the identity outside $\cB_i[\lambda^{-n}, \tau] \subset U_i$.
And
\begin{equation}\label{e.zzz}
H_i \left( \cB_i[a,b] \right) = \cB_i[a, \kappa b] \quad
\text{if $0 < a \le \lambda^{-n} (1-\tfrac{\delta}{2})$ and $b \le (1-\tfrac{\delta}{2})\tau$.}
\end{equation}

It remains to check that~\eqref{e.shrink} holds;
so let $k \le i \le n$.
In the following diagram, $X \xrightarrow{F} Y$ means $F(X) \subset Y$.
Using repeatedly~\eqref{e.lambda}, \eqref{e.chatice}, \eqref{e.sono},
and~\eqref{e.zzz}, 
\begin{multline*}
V_i \subset
\cB_i[\lambda^n (1-\delta) ,(1-\delta) \tau]        \xrightarrow{H_i}
\cB_i[\lambda^n (1-\delta) ,\kappa(1-\delta) \tau] \xrightarrow{L_i} \\
\cB_{i-1}[\lambda^{n+1} (1-\delta) ,\kappa(1-\delta) \tau]
\xrightarrow{L_{i-1} \circ H_{i-1}} \cdots \xrightarrow{L_{i-k+1} \circ H_{i-k+1}} \\
\cB_{i-k}[\lambda^{n+k}(1-\delta) ,\kappa^k (1-\delta)\tau]
\subset
\cB_{i-k}[\lambda^{-n}, \delta \tau] \subset W_{i-k} \, .
\end{multline*}
This proves lemma~\ref{l.slicing}.
\end{proof}

\subsection{Proof of theorem~\ref{t.main}}\label{ss.proof}

Define the following (open) subsets  of $C^1(M,M)$:
\begin{multline*}
\cV_\eps = \big\{f\in C^1(M,M); \;
\text{there exist $K\subset M$ compact, $k\in\N$ such that} \\
\text{$m(K)>1-\eps$  and $m(f^k K)<\eps \big\}$.}
\end{multline*}
By lemma~\ref{l.no acim criterium},
it suffices to show that each $\cV_\eps$ is dense to prove the theorem.
So let $f \in C^1(M,M)$ and $\eps>0$ be fixed;
we will explain how to find $g\in \cV_{4\eps}$ close to $f$.
For clarity we split the proof into steps.

\subsubsection{Step 1: linearizing $f$ on an open tower.}

Let $P_f$ be the set of periodic points of $f$,
and $C_f$ be the set of critical points of $f$.
We can assume (perturbing $f$ if necessary) that $m(P_f \cup C_f)=0$.
(Indeed, it suffices to take $f$ analytic and Kupka-Smale.)

Let $0<\delta<\eps$ be such that $(1-\delta)^d > 1-\eps$. 
Let $k = k(\eps,\delta)$ be given by lemma~\ref{l.slicing}.
Take $n \in \N$ such that $k/(n+1) < \eps$.
Now apply theorem~\ref{t.rokhlin} (and remark~\ref{r.open})
with $\ell = k$, $n_0 = n+1$, $\eps_0 = \eps/2$,
to find an open set $U \subset M$ such that 
$$
\text{$\overline U$, $f^{-1}(\overline U)$, \ldots, $f^{-n}(\overline U)$ are disjoint,}
$$
$$
\sum_{i=0}^{k-1} m(f^{-i} U) < \eps, \qquad
\sum_{i=0}^n m(f^{-i} U) > 1-\eps.
$$

It follows easily from lemma~\ref{l.linearization} that there exist open
sets $Q_i \subset f^{-i}(U)$, $0 \leq i \leq n$ and a $C^1$ perturbation
$\tilde f$ of $f$ such that $\tilde f(Q_i)=Q_{i-1}$, $1 \leq i \leq n$,
$\tilde f|Q_i$ is locally linear and invertible, and $\sum_{i=0}^n
m(Q_i)>1-\epsilon$.  We can assume further (by slight shrinking the $Q_i$)
that each $Q_i$ has only finitely many connected components and $\tilde f$
maps each connected component of $Q_i$ onto a connected component of
$Q_{i-1}$.  We have
\begin{equation}\label{e.rok3}
\sum_{i=0}^{k-1} m(Q_i) < \eps, \qquad
\sum_{i=0}^n m(Q_i) > 1-\eps.
\end{equation}

To simplify writing, we replace $\tilde f$ by $f$.



To simplify things further, we will assume $\bigcup Q_i$ is a subset of $\R^d$,
in order to avoid mentioning the charts.

\subsubsection{Step 2: defining the perturbation $g$.}
Let $k=k(\epsilon,\delta)$ be given by lemma \ref {l.slicing}.
For each sequence $\bar{x} = (x_m, \ldots, x_0)$, $k \leq m \leq n$
with $f(x_i) = x_{i-1}$ and $x_i \in Q_i$,
we apply lemma~\ref{l.slicing} to the sequence of linear maps $Df(x_i)$,
obtaining a certain $\tau_0(\bar x) > 0$.  There are only finitely many
possibilities for the sequence of linear maps $Df(x_i)$, so we can choose
$\tau>0$ such that $\tau<\tau_0(\bar x)$ for all $\bar x$.

For $y \in Q_0$ and small $a$, $b>0$, write
$$
\cB[y,a,b] = y + [-a,a]^{d-1} \times [-b, b] \, .
$$
The family of boxes $\cB[y,r, \tau r] \subset Q_0$ constitute
a Vitali covering
of $Q_0$.
So we can find a finite set $F \subset Q_0$ and
numbers $r(y) > 0$, $y\in F$, such that
$U_0(y) = \cB[y, r(y), \tau r(y)] \subset Q_0$ are disjoint and
\begin{equation} \label {e.vitali}
m\left(Q_0 \setminus \bigsqcup_{y \in F} U_0(y) \right) \leq \epsilon m(Q_0).
\end{equation}
Let $V_0(y)$, $W_0(y)$ be boxes around $y$ as in lemma~\ref{l.slicing}, namely
\begin{align*}
V_0(y) &= \cB[y, (1-\delta) r(y), (1-\delta) \tau r(y)], \\
W_0(y) &= \cB[y, r(y), \delta \tau r(y)].
\end{align*}

Let $\bar F$ be the set of the sequences
$\bar y = (y_m, \ldots, y_0)$ with $k \leq m \leq n$,
$f(y_i)=y_{i-1}$, $y_i \in Q_i$, and $y_0 \in F$.
Then $\bar F$ is finite.
For each $\bar y \in \bar F$ and $i=1,\ldots,m$,
let $U_i(\bar y)$ be the image of $U_0(y_0)$ by the branch of $f^{-i}$
that takes $y_0$ to $y_i$.
Notice that the $U_i(\bar y)$ are either disjoint or coincide, and if
$U_i(\bar y)=U_j(\bar y')$ then $i=j$ and the last $i+1$ symbols in
$\bar y$ and $\bar y'$ coincide.
We define sets $V_i(\bar y)$ and $W_i(\bar y)$ analogously. 

The choice of $\delta$ together with the linearity of $f$ gives
\begin{equation}\label{e.comparison}
\frac{m(V_i(\bar y))}{m(U_i(\bar y))} > 1-\eps \quad\text{and}\quad
\frac{m(W_i(\bar y))}{m(U_i(\bar y))} < \eps \quad\text{for all $\bar{y}$, $i$.}
\end{equation}

Let $H_{i,\bar y}$ be the diffeomorphisms given by lemma~\ref{l.slicing}.
We define $h_{i, \bar y}$ as
$$
h_{i,\bar y}(x) = y_i +
r(y_0) \cdot H_{i, \bar y} \left( (x - y_i )/r(y_0) \right) \, ,
\text{ for $y \in U_i(\bar y)$,}
$$
and $h_{i, \bar y} = \id$ outside $U_i(\bar y)$.  Notice that $h_{i, \bar
y}$ only depends on $U_i(\bar y)$.
Then $h_{i,\bar y}$ is $C^1$-close to the identity.
Moreover, for all $i = k$, $k+1$, \ldots, $m$ we have
$$
f \circ h_{i-k+1, \bar y} \circ \cdots \circ f \circ h_{i-1, \bar y}
\circ f \circ h_{i,\bar y} (V_i(\bar y)) \subset W_{i-k}(\bar y) \, .
$$

We define a perturbation $g:M \to M$ of $f$ as follows:
$g = f \circ h_{i, \bar y}$ on each $U_i(\bar y)$,
and $g$ equals $f$ on
$$
M \setminus \bigcup_{\bar y \in \bar F} \bigcup_{i=0}^n U_i(\bar y) \, .
$$
It follows that
$g^k (V_i(\bar y)) \subset W_{i-k}(\bar y)$
for all $\bar y \in \bar F$.

\subsubsection{Step 3: verifications.}

Define the compact set
$$
K = \bigcup_{\bar y \in \bar F} \bigcup_{i=k}^{n} V_i(\bar y).
$$
First let us see that $K$ has almost full measure.
We have
\begin{multline*}
M \setminus K =
\overbrace{\left(M \setminus \bigcup_{i=0}^n Q_i \right )}^{(\mathrm{I})} \sqcup
\overbrace{\bigcup_{i=0}^n \left(Q_i \setminus
               \bigcup_{\bar y \in \bar F} U_i(\bar y) \right)}^{(\mathrm{II})} \sqcup \\ \sqcup
\underbrace{\bigcup_{i=0}^n \bigcup_{\bar y \in \bar F}
                (U_i(\bar y)\setminus V_i(\bar y))}_{(\mathrm{III})} \sqcup
\underbrace{\bigcup_{i=0}^{k-1} \bigcup_{\bar y \in \bar F} V_i(\bar y)}_{(\mathrm{IV})} \, .
\end{multline*}
From~\eqref{e.rok3} we get $m(\mathrm{I}) < \eps$.
By~\eqref{e.vitali} and linearity of $f$, we have $m(\mathrm{II})< \eps$.
From~\eqref{e.comparison}, $m(\mathrm{III}) \le \eps$.
Finally, using~\eqref{e.rok3},
$$
m(\mathrm{IV})
\le m\left(\bigcup_{i=0}^{k-1}  Q_i \right ) < \eps.
$$
So we obtain that $m(M \setminus K) < 4\eps$.

Next let us see that $g^k K$ has small measure.
We have
$$
g^k K = \bigcup_{\bar y \in \bar F} \bigcup_{i \geq k} g^k V_i(\bar y)
\subset \bigcup_{\bar y \in \bar F} \bigcup_{i \geq 0} W_i(\bar y).
$$
Using~\eqref{e.comparison},
$$
m(g^k K)
\le \eps m \left( \bigcup_{\bar y \in \bar F} \bigcup_{i \geq 0} U_i(\bar y) \right)
< \eps.
$$

We have shown that $g \in \cV_{4\eps}$.
This proves theorem~\ref{t.main}.

\begin{ack}
We thank Bassam Fayad and Enrique Pujals
for posing to us the problem addressed in this paper.
The results were obtained during a visit of the second author to the Coll\`ege de France,
which was supported by the Franco-Brazilian cooperation agreement in Mathematics.
\end{ack}

\end{document}